\documentclass[12pt]{article}

\usepackage{amsmath,epsfig,epsf,psfrag,natbib}
\usepackage{amssymb}
\usepackage[latin1]{inputenc}
\usepackage{rotating}
\usepackage{subfigure}
\usepackage{hyperref}
\usepackage{amsbsy}
\usepackage[ruled,vlined]{algorithm2e}
\usepackage{lscape,url, epstopdf}
\usepackage{color,latexsym,amsfonts,amsthm,amscd,graphicx}
\usepackage{mathrsfs}
\usepackage{float, multirow, booktabs, threeparttable, diagbox}
\usepackage{algorithmic}
\newcommand{\ben}{\begin{enumerate}}
	\newcommand{\een}{\end{enumerate}}
\newcommand{\be}{\begin{equation}}
	\newcommand{\ee}{\end{equation}}
\newcommand{\bas}{\begin{eqnarray*}}
	\newcommand{\eas}{\end{eqnarray*}}
\newcommand{\ba}{\begin{eqnarray}}
	\newcommand{\ea}{\end{eqnarray}}
\newcommand{\bit}{\begin{itemize}}
	\newcommand{\eit}{\end{itemize}}

\newtheorem{theorem}{Theorem}
\newtheorem{lemma}{Lemma}

\newtheorem{remark}{Remark}

\newtheorem{assumption}{Assumption}

\newtheorem{prop}{Proposition}

\newcommand{\e}{ { \mathbb{E}}}

 \def\T{{ \mathrm{\scriptscriptstyle \top} }}

\newcommand{\sha}{{\rm sha}}
\newcommand{\pr}{ {\rm pr} }

\newcommand{\convergeto}{ {\overset{d}{\longrightarrow \; }}}

\newcommand{\Perp}{\perp\kern-7pt\perp}
\newcommand{\notPerp}{ \not  \kern-2pt \perp\kern-7pt\perp}

\usepackage{authblk}

\begin{document}
\date{}
\title{
Efficient estimation of partially linear additive Cox models and variance estimation under shape restrictions
}
\author[1]{Junjun Lang}
\author[1]{Yukun Liu\thanks{Corresponding author:  ykliu@sfs.ecnu.edu.cn}}
\author[2]{Jing Qin}
{\small   \affil[1]{ \small
		KLATASDS-MOE,  School of Statistics,
		East China Normal University,
		Shanghai 200062, China}
	\affil[2]{ \small
		National Institute of Allergy and Infectious Diseases, National Institutes of Health, USA}
	}

\maketitle
\abstract{

Shape-restricted inferences  have
exhibited empirical success in various applications with survival data.
However, certain works fall short in providing a rigorous theoretical justification
and an easy-to-use variance estimator with theoretical guarantee.
Motivated by \cite{deng2023active},  this paper
delves into an additive and shape-restricted partially linear Cox model for right-censored data,
where each additive component satisfies a specific shape restriction,
encompassing monotonic increasing/decreasing and convexity/concavity.
We systematically investigate the   consistencies and convergence rates
of the shape-restricted maximum partial likelihood estimator (SMPLE)
of all the underlying parameters.
We further establish the aymptotic normality  and semiparametric effiency
of  the SMPLE for the linear covariate shift.
To estimate the asymptotic variance, we propose an innovative data-splitting
variance estimation method that boasts exceptional versatility and broad applicability.
  Our simulation results and an analysis of
the Rotterdam Breast Cancer   dataset demonstrate that
the SMPLE  has comparable performance with the maximum likelihood estimator
under the Cox model  when the Cox model is correct,
and outperforms the latter  and \cite{Huang1999}'s method
 when the Cox model is violated or the hazard is nonsmooth.
Meanwhile, the proposed variance estimation method
usually leads to  reliable interval estimates
based on  the SMPLE and its competitors.

	
}

\section{Introduction}

Shape restrictions (such as monotonicity and convexity) arise naturally in numerous practical scenarios.
For instance, the growth curves of animals and plants in ecology and
 the dose-response in medicine must inherently
exhibit non-decreasing characteristics \citep{chang2007shape, wang2012shape}.
In the realm of economics, utility and production functions are often
concave in income and prices  \citep{matzkin1991semiparametric, varian1984nonparametric},
cost functions are monotone increasing, concave in input prices, and may exhibit non-increasing or
non-decreasing returns to scale  \citep{horowitz2017nonparametric}.
In genetic epidemiology studies, the cumulative risk of a disease
for individuals  possess monotonicity \citep{qin2014combining}.
 While in reliability analysis, the bathtub curve describing the failure rate typically displays
convexity.

Incorporating shape restrictions into statistical analysis,
apart from its exceptional interpretability and ability to
enforce domain-specific constraints, often results in an estimation procedure
that is devoid of tuning parameters, enhancing its efficiency
 and robustness. Therefore shape-restricted techniques has become
 an increasing popular tool for statistical inference or learning
 in various settings over the past decades.
  A comprehensive review on shape-restricted nonparametric inferences
 can be found in  \cite{groeneboom2014} and references therein.
Recently, \cite{ChenandSamworth2016} developed an algorithm for
 the estimation of the generalized additive model
 in which each of the additive components is linear or subject to
 a shape restriction.
 \cite{balabdaoui2019score} considered the estimation of
 the index parameter  in a single-index model with a monotonically increasing link function.
 \cite{DengandZhang2020} studied   minimax and adaptation rates in general multiple isotonic
regression.
 \cite{feng2022nonparametric} systematically investigate
 the theoretical properties of the least squared estimator
  of a S-shaped regression function.

This paper focus on the statistical inference for righter-censored
survival data. Let $T$ denote  the survival time and $(Z,X)\in\mathbb{R}^{p}\times\mathbb{R}^{d}$
denote  a $(p+d)\times 1$ vector of covariates.
We consider the partially linear Cox model of \cite[PLCM]{Sasieni1992}
for modelling the conditional hazard function, i.e.
\ba
\label{model-cox-additive}
\lambda_{T}(t\mid x, z) = \lambda(t)\exp(\beta^{\top}x+g(z)),
\ea
where $\lambda(\cdot)$ is the unspecified baseline hazard function, $\beta\in\mathbb{R}^{d}$ is unspecified and $g(\cdot):\mathbb{R}^{p}\mapsto\mathbb{R}$ is an unknown function.
This model   reduces to
the renowned Cox proportional hazards model \citep{Cox1972, Cox1975}
when   the covariate $Z$ disappears,
and it  becomes the nonparametric Cox
model \citep{sleeper1990regression, o1993nonparametric}
in the absence of $X$.

Many nonparametric techniques have been developed
for the estimation of the PLCM, in particular for the linear covariate effect.
Examples include   profile partial likelihood
together with a kernel technique
\citep{heller2001cox},
maximum likelihood estimation with
a deep neural network
\citep{Zhong2022},
and  a kernel machine representation method
\citep{rong2024kernel}, etc.
However, these methods are either hampered
by the curse of dimensionality or lack interpretability
for  $g(Z)$,  or suffer from tuning parameters,
whose selection is not always straightforward.
Alternatively,  \cite{Huang1999} proposed to model  $g(Z)$ by a generalized
additive model \citep{1986Generalized, 1990Generalized},
which effectively  avoids the curse of dimensionality
and enforces an additive effect for the covariate $Z$. Specifically,
\ba
\label{model-g}
g(Z) =\sum\limits_{j=1}^{p}g_{j}(Z_{(j)}),
\ea
where  for $1\leq j\leq p$,
$Z_{(j)}$ is the $j$-th component of $Z$ and $g_{j}$ is an unknown function.
\cite{Huang1999} proposed the use of polynomial splines to fit the unknown additive components.
This method   entails a   number of tuning parameters,
  also yields convergence rates that lack conciseness and elegance.
Furthermore,  the spline method  does not provide good interpretability
for the additive covariate $Z_{(j)}$.

Our paper is motivated by the work of
\cite{deng2023active}
which studied  a shape-restricted  and  additive PLCM.
Specifically, under  models  \eqref{model-cox-additive}  and \ref{model-g},
 they assume that   each $g_{j}$ is
monotonic  increasing/decreasing or convex/concave.
An active-set optimization algorithm was provided to
calculate the shape-restricted maximum likelihood estimator.
The shape-restriction strategy facilitates the utilization of prior knowledge
regarding the effect of the log conditional hazard function
on each covariate $Z_{(j)}$ and leads to a  tuning-parameter-free estimation procedure.
However, they  proved only a consistency result, and  did not provide any asymptotic
normality results.
\cite{qin2021estrogen} studied a PLCM
with a single additive component subject  to shape restrictions,
but   they did not establish any $\sqrt{n}$-consistency result.
In addition, in general shape-restriction inferences,
even if asymptotic normality results can be established,
it is generally challenging to construct
reasonable estimators for the asymptotic variances
with theoretical guarantees \citep{Groeneboom2017}.

This paper makes two main contributions
to the literature of  additive and shape-restricted PLCMs for survival data.
The first contribution   is to
provide powerful statistical guarantees for the  shape-restricted maximum
partial likelihood estimator (SMPLE) and
the induced Breslow-type estimator for the baseline cumulative hazard function
under the model assumption of \cite{deng2023active}.
This includes a thorough convergence rate analysis for the
estimators of the infinitely dimensional parameters,
as well as establishing asymptotic normality and semiparametric
efficiency for the estimator of the linear covariate effect.
Our second contribution is to  offer an easy-to-use  estimator for
the asymptotic variance of the linear covariate effect estimator.
We show that this variance estimation method always
provide consistent estimators once
the corresponding asymptotic normality result holds.
This method is very flexible and is applicable for general purpose
especially in  shape-restricted inferences,
where theoretical guarantee of a bootstrap variance estimator
is generally rather challenging  \citep{Groeneboom2017}.
Our simulation results and an analysis of
the Rotterdam Breast Cancer   dataset demonstrate that
the SMPLE  has comparable performance with the maximum likelihood estimator
under the Cox model  when the Cox model is correct,
and outperforms the latter  and \cite{Huang1999}'s method
 when the Cox model is violated and the hazard is nonsmooth.
Meanwhile,
the proposed variance estimation method
usually leads to  reliable interval estimates
for the SMPLE and its competitors.


The rest of this paper is organized as follows.
Section \ref{sec-model} introduces notations,
data, and  the shape-restricted maximum partial
likelihood estimators (SMPLE).
Section \ref{sec-rate-convergence}
investigates   the convergence rates of
 the SMPLEs for all the unknown parameters,
including $\beta$,
the unknown additive components,
and the baseline cumulative hazard function.
Section
  \ref{sec-asym-nor-effi} establishes
the asymptotic normality and semiparametric efficiency of
the SMPLE for  $\beta$.
A novel  estimation method is also provided
to estimate the asymptotic variance of the SMPLE of $\beta$.
 A simulation study
    and real data analysis are presented
    in Section \ref{sec-sim} and \ref{sec-real-data},
     respectively.
 Section \ref{sec-con} contains  concluding remarks.
 For clarity, all technical proofs are postponed to the
 supplementary material.

\section{Methodology}
\label{sec-model}
\subsection{Data and model assumptions}

Let   $T $ and $(X, Z)$  be the survival time  and  the vector of covariates, respectively,
 in the introduction.
Suppose that  given  $(X, Z)$,
the conditional hazard function of   $T$
satisfies model \eqref{model-cox-additive} with $g(Z)$
satisfying \eqref{model-g}.
The survival time  $T $  may be right censored
by a censoring time $C$
and we only observe $Y=\min(T, C)$.
Throughout this paper,
we use $\mathbf{1}(A)$ to denote the indicator function of the set $A$
and use a subscript 0 to highlight the true counterpart of a parameter.
Let $\Delta = \mathbf{1}(T\leq C)$ be the non-censoring indicator.
Given $n$ independent and identically distributed (iid)
observations  $(X_i, Z_i, Y_i, \Delta_i)$,  $1\leq i\leq n$,
from $(X, Z, Y, \Delta)$,
we wish to infer $(\beta_{0}, g_{0})$ and the baseline cumulative hazard function $\Lambda_{0}(y) = \int_{0}^{y}\lambda_{0}(t) dt$.

The identifiability issue of models \eqref{model-cox-additive} and \eqref{model-g}
was investigated by \cite{deng2023active}, following which
  we assume  $\e\{ g_{0,j}(Z_{(j)})  \Delta\}=0$, $j=1, 2, \cdots, p$,
for identifiability.
Furthermore, we assume that for  $1\leq j \leq p$,
$g_{0,j}(\cdot)$ satisfies one of the four shape restrictions:
monotone increasing, monotone decreasing,
 convex and concave, which are encoded as shape types 1, 2, 3 and 4, respectively.
For any additive function $g=\sum\limits_{j=1}^{p}g_{j}(z_{(j)})$, we define $\sha(g)=(\sha(g_{1}), \cdots, \sha(g_{p}))^{\top}$, where $\sha(h)\in\{1, 2, 3, 4\}$ denotes the shape type of a univariate function $h$.
We always denote $\boldsymbol{k}_{0}=\sha(g_{0})$.
Let $\mathcal{X}$ be the support of $X$   and
for simplicity, we assume that the support of $Z_{(j)}$ is $[0, 1]$ for $j=1,2, \cdots,p$.

\subsection{SMPLE}

For any $(\beta, g)$, denote $\eta = (\beta, g)$ and
$R_{\eta}(U) = X^{\top}\beta+g(Z)$, where $U = (X^{\top}, Z^{\top})^{\top}$.
$1/n$ times the usual partial log likelihood
is
\bas
L_{n}(\eta) = \frac{1}{n}\sum\limits_{i=1}^{n}\Delta_{i}\left[R_{\eta}(U_{i})-\log\left(\sum\limits_{j=1}^{n}\mathbf{1}(Y_{j}\geq Y_{i})\exp(R_{\eta}(U_{j}))\right)\right].
\eas
We propose to estimate
$\eta$ by
the shape-restricted maximum partial likelihood estimator  (SMPLE),
\ba
\label{definition-estimator}
\hat{\eta}:=(\hat{\beta}, \hat{g}) = \operatorname{argmax}_{\eta\in\mathbb{R}^{d}\times\mathcal{G}_{\boldsymbol{k}_{0}}}L_{n}(\eta),
\ea
where
\begin{equation*}
	\begin{aligned}
	\mathcal{G}_{\boldsymbol{k}_{0}} = \bigg\{g:[0, 1]^{p}\mapsto\mathbb{R} \mid g(Z) = \sum\limits_{j=1}^{p}g_{j}(z_{(j)}),~ \sha(g) = \boldsymbol{k}_{0}, \mathbb{E}\left[\Delta g_{j}(Z_{(j)})\right]=0,~1\leq j\leq p\bigg\} 	
	\end{aligned}
\end{equation*}
is the parameter space of $g$.
With the SMPLE in \eqref{definition-estimator}, we   estimate $\Lambda_{0}(y)$  by the Breslow-type estimator
\ba
\label{estimator-cumulative-hazard-function}
\hat{\Lambda}(y; \hat{\eta}) = \frac{1}{n}\sum\limits_{j=1}^{n} \frac{\Delta_{j}}{ S_{0,n}(Y_{j}, \hat{\eta})} \mathbf{1}(y\geq Y_{j}) ,
\ea
where
$
S_{0,n}(y, \eta) =  (1/n)\sum\limits_{i=1}^{n}\{  \mathbf{1}(Y_{i}\geq y)\exp(R_{\eta}(U_{i}))\}.
$

The SMPLE defined in \eqref{definition-estimator} can be calculated with the active-set algorithm introduced in \cite{deng2023active}.
 Let  $\hat{g}_{j}$ be functions satisfying   $\hat{g}(Z) = \sum\limits_{j=1}^{p}\hat{g}_{j}(Z_{(j)})$
for all $Z$.
The function $\hat{g}(Z)$ is  unique only at the observed $Z_i$
and  is therefore   non-unique typically for $Z$ other than $Z_i$'s,
 which is akin to  general shape-restricted regression estimators
 \citep{ChenandSamworth2016}.
This implies that   $\hat{g}_{j}(\cdot)$ is usually  non-unique for $1\leq j\leq p$,
and the solution set of $\hat{g}_{j}(\cdot)$ always contains a piece-wise linear function
 \citep{deng2023active}. See Figure \ref{real-plot} for an illustration.

\section{Rate of convergence}
\label{sec-rate-convergence}

The consistency property of the SMPLEs  $\hat{\beta}$,
$\hat{g}$, and $\hat{g}_{j}$ was established by \cite{deng2023active}.
In this section, we establish their   convergence  rates.
We make the following assumptions.

\begin{assumption}
\label{assum-basic}
(i) The observed $\{Y_{i}\}_{i=1}^{n}$ are in the interval $[0,\tau]$, for some $\tau>0$.
(ii) Given $U$, $T$ and $C$ are mutually independent of each other.
(iii) $\Lambda_{0}(\tau)<\infty$ and $\pr(C\geq \tau\mid U)\geq c>0$ almost surely  for some constant $c$.
(iv) $\mathbb{E}[\Delta X]=0$ and $\mathbb{E}[\Delta]>0$.
\end{assumption}

Let $\|\cdot\|$ denote the usual Euclidean norm and $\|f(\cdot)\|_{\infty}$
 the supreme norm of a real-valued function $f$.
For any constant $M>0$,  define
\bas
\mathcal{K}_{M, \boldsymbol{k}_{0}}:=\left\{\eta\mid \eta = (\beta, g), g\in\mathcal{G}_{\boldsymbol{k}_{0}}, \|\beta\|+\sum\limits_{j=1}^{p}\|g_{j}(\cdot)\|_{\infty}\leq M\right\}.
\eas
\begin{assumption}
	\label{assum-bounded}
The support 	$\mathcal{X}$ of $X$ is a bounded subset of $\mathbb{R}^{d}$ and there exists a positive constant $M_{0}>0$  such that $\eta_{0}\in\mathcal{K}_{M_{0}, \boldsymbol{k}_{0}}$.
\end{assumption}
\begin{assumption}
	\label{assum-censor prob}
	There exists a small positive constant $\epsilon$ such that $\pr(\Delta = 1\mid U)>\epsilon$ almost surely with respect to the probability measure of $U$.
\end{assumption}

\begin{assumption}
	\label{assum-B5-Huang}
	The joint density  of $(Y, Z, \Delta )$ satisfies
	\bas
	0<\inf_{(y, z)\in [0,\tau]\times[0, 1]^{p}}\pr(Y=y, Z=z, \Delta=1)\leq \sup_{(y, z)\in [0,\tau]\times[0, 1]^{p}}\pr(Y=y, Z=z, \Delta=1)<\infty.
	\eas
\end{assumption}

\begin{assumption}
	\label{assum-convex-bounded-density}
	When ${\sha}(g_{0,j})\in\{3,4\}$, the density function $Z_{(j)}$ with respect to
 the Lebesgue measure has uniformly upper and lower bounds on $[0, 1]$.
\end{assumption}

Assumptions \ref{assum-basic}--\ref{assum-bounded} are  standard in the theoretical analysis of traditional Cox model and its variants \citep{Huang1999, Zhong2022}.
Assumption \ref{assum-censor prob} ensures that the probability of being uncensored is positive regardless of the covariate values, and it  is used to establish  the convergence rate results
in Theorem \ref{thm-rate-convergence-SMPLE}.
Assumption  \ref{assum-B5-Huang} is used in the calculation of
the semiparametric efficiency lower bound  \citep{Huang1999}.
In Assumption \ref{assum-convex-bounded-density},
the upper bound requirement  is used to calculate some entropy results
needed in the proof of Proposition \ref{prob-consistency},
and the lower bound requirement  guarantees
that the approximation errors of
piecewise linear approximations of
the convex/concave additive components
to themselves are small enough   in the proof of Theorem \ref{thm-asym-normal-linear}.

The Proposition below establishes the consistency of $R_{\hat{\eta}}(\cdot)$,
as an estimator of $R_{\eta_{0}}(\cdot)$,
which roughly implies the consistency of $\hat{\eta}=(\hat{\beta},\hat{g})$.
We would have  proved the  consistencies of
$\hat{\beta}$ and each $\hat{g}_{j}$ separately.
However, the latter results are not needed in the proofs of the subsequent
convergence rate results given  the consistency of $R_{\hat{\eta}}(\cdot)$.

\begin{prop}
	\label{prob-consistency}
Suppose that  models \eqref{model-cox-additive} and \eqref{model-g}
 and Assumptions \ref{assum-basic}, \ref{assum-bounded} and \ref{assum-convex-bounded-density} are satisfied.
 As $n\rightarrow \infty$, we have
	\ba
	\label{consistency-together}
	\|R_{\hat{\eta}}(\cdot)-R_{\eta_{0}}(\cdot)\|_{\infty}=o_{p}(1).
	\ea
\end{prop}

Define
$
d^{2}(\eta, \eta_{0}) = \mathbb{E}_{U}\left\{(R_{\eta}(U)-R_{\eta_{0}}(U))^{2}\right\},
$
where $\mathbb{E}_{U}$ denotes the expectation with respect to $U$.
Let $\|\cdot\|_{L_{2}}$ denote the $L_{2}(P)$ norm and $\rho(\boldsymbol{k}_{0}) = 0.5+0.5 \cdot \mathbf{1}( \cup_{i=1}^p (\sha(g_{0,i})\in\{1,2\} ) )$.
One of our main results is to establish  the  convergence rate of the SMPLE $\hat{\eta}$.
\begin{theorem}
			\label{thm-rate-convergence-SMPLE}
Assume the same conditions in Proposition \ref{prob-consistency}.  As $n\rightarrow \infty$, we have
\[
	d(\hat{\eta},\eta_{0}) = O_{p} \left(n^{-\frac{1}{2+\rho(\boldsymbol{k}_{0})}} \right).	
\]
Furthermore, if Assumptions \ref{assum-censor prob}--\ref{assum-B5-Huang} are satisfied and
$I(\beta_{0})$ (defined in \eqref{def-Ibeta0}) is non-singular,
then   for $1\leq j\leq p$,
	\bas
	\|\hat{\beta}-\beta_{0}\| = O_{p}\left(n^{-\frac{1}{2+\rho(\boldsymbol{k}_{0})}}\right),\quad \|\hat{g}_{j}(Z_{(j)})-g_{0,j}(Z_{(j)})\|_{L_{2}} = O_{p}\left(n^{-\frac{1}{2+\rho(\boldsymbol{k}_{0})}}\right).
	\eas	
\end{theorem}

According to Theorem \ref{thm-rate-convergence-SMPLE},
the rates of convergence of all $\hat{g}_{j}$ are $O_{p}(n^{-2/5})$
if none of the additive components of $g$ is   monotonic.
Conversely, if one  additive component of $g$ is  monotonic,
then their convergence rates all slow down to $O_{p}(n^{-1/3})$.
An explanation for this finding is that the complexity of the class of bounded
and monotonic functions is much larger than that of the class of bounded and convex
(or concave) functions. These convergence rate results are free from
the covariate dimensionality and exhibit a much more elegant form
than those in \cite{Huang1999} and \cite{Zhong2022}.
Theorem \ref{thm-rate-convergence-SMPLE} also establishes
the convergence rate of $\hat{\beta}$, although it is sub-optimal.

With Theorem \ref{thm-rate-convergence-SMPLE}, we are able to
establish the uniformly rate of convergence
for the SMPLE $\hat{\Lambda}(y; \hat{\eta})$  in \eqref{estimator-cumulative-hazard-function}
of the baseline cumulative hazard function $\Lambda_{0}(y)$.
It turns out that   $\hat{\Lambda}(y; \hat{\eta})$
has the same convergence rate as  $\hat \eta$, although
their convergence rates are quantified by different distances.

\begin{theorem}
	\label{thm-uniformly-rate-hazard}
Assume the same conditions as in Proposition \ref{prob-consistency}.
As $n\rightarrow \infty$, it holds that
	\bas
	\sup_{y\in[0,\tau]}\left|\hat{\Lambda}(y; \hat{\eta})-\Lambda_{0}(y)\right| = O_{p}\left(n^{-\frac{1}{2+\rho(\boldsymbol{k}_{0})}}\right).
	\eas	
\end{theorem}

\begin{remark}
	\label{remark-convex-monotone}
In practice,  one may impose a combination of monotonicity and convexity/concavity constraints
on the additive components according to domain knowledge.
 See \citep{ChenandSamworth2016, Kuchibhotla2023, deng2023active} for further motivation
on additional shape constraints.   Proposition \ref{prob-consistency}
 and Theorems \ref{thm-rate-convergence-SMPLE}--\ref{thm-variance-consistency} still hold
when model \eqref{model-g} incorporates additive components
that satisfy both monotonicity and convexity/concavity restrictions.
An intuitive explanation for this result is that the parameter space 	$\mathcal{G}_{\boldsymbol{k}_{0}}$ is
reduced by  additional constraints on the additive components
and this can lead to better convergence rates of the SMPLE (if not the same).
\end{remark}

\section{Asymptotic normality and efficiency}
\label{sec-asym-nor-effi}

Based on the convergence rate results in the previous section,
in this section, we further  show that our SMPLE  in \eqref{definition-estimator}
for   the linear covariate effect  $\hat{\beta}$
is asymptotically normal and semiparametric efficient,
in the sense that its asymptotic variance achieves
the semiparametric efficiency lower bound   \citep{bickel1993efficient}
or the information bound of
  estimating $\beta_{0}$ under models \eqref{model-cox-additive} and  \eqref{model-g}.

We begin with presenting  the information bound for $\beta_{0}$. Recall that $U = (X^{\top}, Z^{\top})^{\top}$,
and  define
\bas
M(y)\equiv M(y\mid Y, \Delta, U) = \Delta \mathbf{1}(Y\leq y) - \int_{0}^{y}\mathbf{1}(Y\geq t)\exp\{ R_{\eta_{0}}(U) \} d\Lambda_{0}(t),
\eas
which is a  counting process martingale associated with the Cox model.
The log-likelihood of model \eqref{model-cox-additive}
based on one observation $( X, Z,Y, \Delta)$ is  (up to constant)
\ba
\label{log-likelihood}
\ell(\beta, g, \Lambda) = \Delta\log\lambda(Y) + \Delta\{X^{\top}\beta+g(Z)\}-\Lambda(Y)\exp\{X^{\top}\beta+g(Z)\}.
\ea
Conisder a parametric smooth sub-model $\{\lambda_{(\nu)}:\nu\in\mathbb{R}\}$ and $\{g_{j, (\nu)}: \nu \in\mathbb{R}\}$, $1\leq j\leq p$, with  $\lambda_{(0)} = \lambda_{0}$ and $g_{j, (0)} = g_{0,j}$.
Define $L_{2}(P_{Y})$ to be the set of   $a(\cdot)$ satisfying
 $\e \{ \Delta a^{2}(Y)\}<\infty$
 and
$
 a(y) =  \partial\log \lambda_{(\nu)}(y)/\partial\nu  |_{\nu=0},
$
  for some submodel.
Similarly, for $ 1\leq j \leq p$,
define  $L_{2}^{0}(P_{Z_{(j)}})$
to be the set of   $h_{j}$ satisfying $\e\{ \Delta h_{j}(Z_{(j)}) \} = 0$,
$\e\{ \Delta h_{j}^{2}(Z_{(j)}) \} <\infty$ and
$
h_{j}(z_{(j)})=
 \partial g_{j, (\nu )}(z_{(j)}) / \partial \nu  |_{\nu =0}
$
 for some submodel.
The following lemma, which is  Theorem 3.1 of \cite{Huang1999},
gives the information bound of $\beta_0$.

\begin{lemma}[Theorem 3.1 of \cite{Huang1999}]
	\label{lemma-Huang}
Suppose that  models \eqref{model-cox-additive} and \eqref{model-g} and
Assumptions \ref{assum-basic}--\ref{assum-B5-Huang} are satisfied.
Let
$((\boldsymbol{a}^{\star})^{\top}, (\boldsymbol{h}_{1}^{\star})^{\top}, \cdots, (\boldsymbol{h}_{p}^{\star})^{\top})^{\top}$
be the unique, vector-valued function  in $L_{2}(P_{Y})^{d}\times L_{2}^{0}(P_{Z_{(1)}})^{d}\times\cdots\times L_{2}^{0}(P_{Z_{(p)}})^{d}$ that minimizes
\bas
\e\left\{
\Delta\|X-\boldsymbol{a}(Y)-\boldsymbol{h}_{1}(Z_{(1)})-\cdots-\boldsymbol{h}_{p}(Z_{(p)})\|^{2}\right\}.
\eas
\bit
\item[(1)]
The efficient score for estimation of $\beta_{0}$  is
	\bas
	\ell_{\beta_{0}}^{\star}(Y, \Delta, U) = \int_{0}^{\tau}\{ X-\boldsymbol{a}^{\star}(y)-\boldsymbol{h}^{\star}(Z) \} dM(y),
	\eas
	where  $\boldsymbol{h}^{\star}(Z) = \sum\limits_{j=1}^{p}\boldsymbol{h}_{j}^{\star}(Z_{(j)})$ and
$
	\boldsymbol{a}^{\star}(y) = \e\left\{
X-\boldsymbol{h}^{\star}(Z)\mid Y=y, \Delta =1\right\}.
$
\item[(2)] The information bound for estimation of $\beta_{0}$ is
	\ba
	\label{def-Ibeta0}
	I(\beta_{0}) = \e\left[
 \left\{ 	\ell_{\beta_{0}}^{\star}(Y, \Delta, U)\right\}^{\otimes 2} \right]
 = \e\left[
 \Delta \left\{
 X- \boldsymbol{a}^{\star}(Y) - \boldsymbol{h}^{\star}(Z) \right\}^{\otimes 2}\right],
 	\ea
where  $A^{\otimes 2} = AA^{\top}$ for any vector or matrix $A$.
\eit
\end{lemma}

Additional assumptions are needed
to obtain the asymptotic normality and efficiency of $\hat \beta$.
Denote $\boldsymbol{h}_{j}^{\star}  =(\boldsymbol{h}_{j, 1}^{\star}, \cdots, \boldsymbol{h}_{j, d}^{\star})^{\top}$  for $1\leq j\leq p$.

\begin{assumption}
	\label{assum-mono}
	When $\sha(g_{0,j})\in\{1,2\}$, there exist  constant $\tilde{C}_{1}>0$ and $\tilde{C}_{2}>0$  such that
$
	\|\boldsymbol{h}^{\star}_{j}(x_{1}) - \boldsymbol{h}^{\star}_{j}(y_{1})\|\leq \tilde{C}_{1}|x_{1}-y_{1}|
$
and
$|g_{0,j}^{-1}(x_{2})-g_{0,j}^{-1}(y_{2})|\leq \tilde{C}_{2}|x_{2}-y_{2}|.
$
for all $x_{1}, y_{1}\in[0, 1]$ and $x_{2}, y_{2}\in[-M_{0}, M_{0}]$.
\end{assumption}

\begin{assumption}
	\label{assum-convex-holder-continous}
	When $\sha(g_{0,j})\in\{3,4\}$,  the function $\boldsymbol{h}_{j, \tilde{i}}^{\star}(\cdot)$ is $\bar{p}$-H$\rm \ddot{o}$lder continuous
for  all $1\leq \tilde{i}\leq d$   and some  $\bar p\in (\rho(\boldsymbol{k}_{0}), \;  2]$.
\end{assumption}

\begin{assumption}
	\label{assum-convex-strict-convex}
	When $\sha(g_{0,j})\in\{3\}$, the function $g_{0,j}(t)-\tilde{C}t^{2}$ is convex for some constant $\tilde{C}>0$;   when $\sha(g_{0,j})\in\{4\}$, the function $g_{0,j}(t)+\tilde{C}t^{2}$ is concave for some constant $\tilde{C}>0$.
\end{assumption}

Assumption \ref{assum-mono} is used to control
 the fluctuation of  the score function
  corresponding to the monotone components
in the direction of the projection defined in Lemma \ref{lemma-Huang}.
 \citep{Huang2002,Cheng2009} adopted a similar assumption.
Assumptions \ref{assum-convex-holder-continous}--\ref{assum-convex-strict-convex},
which are analogues of  Assumption B1 of \cite{Kuchibhotla2023},
are used for approximations of the convex (concave) components.


\begin{theorem}
	\label{thm-asym-normal-linear}
Suppose that  models \eqref{model-cox-additive} and \eqref{model-g}
and Assumptions \ref{assum-basic}--\ref{assum-convex-strict-convex}
are satisfied. Further assume that $I(\beta_{0})$ is non-singular.
Then as $n\rightarrow \infty$,
$
 \sqrt{n} (\hat{\beta}-\beta_{0} )
\stackrel{d}{\longrightarrow} N(0, I^{-1}(\beta_{0})).
$
This   implies that  the asymptotic variance
achieves the information  bound and $\hat{\beta}$
is asymptotically semiparametric efficient among
all regular estimators of $\beta_{0}$.
\end{theorem}

By  Theorem \ref{thm-asym-normal-linear},
$\hat{\beta}$ has an asymptotically normal distribution with asymptotic variance $I^{-1}(\beta_{0})$.
When  making inference about $\beta_0$ based on this theorem,
we need to construct a consistent estimator for $I^{-1}(\beta_{0})$.
However, $I^{-1}(\beta_{0})$ or equivalently $I (\beta_{0})$
has   a rather complicated form, making  its plug-in estimator not easy to use.
To crack this nut,  we propose a  novel  data-splitting
estimation method to estimate $I^{-1}(\beta_{0})$.

  \subsection{Data-splitting  variance estimation and inference on $\beta$}

We introduce the proposed data-splitting variance
estimation method  under a general setting
as it is applicable generally.
Let  $\theta_{0}\in\mathbb{R}^{d}$ be a functional of
 a  statistical population $\mathscr{P}$
 and
 $\hat{\theta}_{n}$ be   an estimator of $\theta_{0}$
based on i.i.d samples $\{O_{i}\}_{i=1}^{n}$ from $\mathscr{P}$.
Suppose that $\sqrt{n} (\hat \theta_n - \theta_{0}  ) \convergeto N(0, \Sigma)$,
where $\Sigma$ is a semi-positive matrix.
Let $k_{n}<n$ and $k_{n}\rightarrow\infty$.
We partition the sample  into $\lfloor k_{n}\rfloor$  subsamples,  each of which
has $m_n = \lfloor n/k_n \rfloor $ observations,
and  let  $\hat \theta_{ni} $ denote the  estimator of $\theta_{0}$ based on
the $i$-th subsample,  $1\leq i\leq \lfloor k_n\rfloor$.
Our splitting-data estimator for the asymptotic variance $\Sigma$ is defined as
\ba
\label{hat-Sigma}
\hat \Sigma  =  \frac{m_n}{\lfloor k_{n}\rfloor}
\sum_{i=1}^{\lfloor k_{n}\rfloor}
( \hat \theta_{ni}  -  \bar{\hat{\theta}}_{n}   )( \hat \theta_{ni}  -  \bar{\hat{\theta}}_{n}   )^{\top},
\ea
where $\bar{\hat{\theta}}_{n}$ is the sample mean of
$\hat \theta_{n1}, \ldots, \hat \theta_{n\lfloor k_n\rfloor}$.
For better stableness,  we may repeat the above splitting and estimating procedure
for many times  and take the average of the resulting variance estimates as
a final variance estimate.

\begin{theorem}
	\label{thm-variance-consistency}
Let  $\hat{\theta}_{n}$ be   an estimator of $\theta_{0}$
based on i.i.d samples $\{O_{i}\}_{i=1}^{n}$ and $\sqrt{n} (\hat \theta_n - \theta_{0}  ) \convergeto N(0, \Sigma)$.
	Let $k_n = n^{\tilde{\alpha}}$ for some $\tilde{\alpha} \in (0, 1)$,
$m_n = \lfloor n^{1-\tilde{\alpha}}\rfloor$ and $\hat \Sigma$
be the variance estimator in \eqref{hat-Sigma}.
Then    $\hat \Sigma  = \Sigma + o_p(1)$ as $n\rightarrow \infty$.
	
\end{theorem}

Theorem \ref{thm-variance-consistency} guarantees the validity of the
data-splitting  variance estimator.
This method is very easy to use
and  is    flexible enough for general purpose.
Alternatively, we may construct a variance estimator by bootstrap.
However, the consistency of a bootstrap variance estimator
often requires stronger conditions \citep{Groeneboom2017}
and is often very difficult to prove, especially
under shape restrictions.

As a specific application,
we apply  the data-splitting  estimation method
to construct an  estimator for
the information bound or the asymptotic variance  $ I^{-1}(\beta_0) $  of $\hat \beta$.
Denote the resulting  estimator by  $\widehat{I^{-1}(\beta_0)}$,
which is consistent to   $ I^{-1}(\beta_0) $  by  Theorem \ref{thm-variance-consistency}.
Therefore  $n  (\hat \beta - \beta_0)^\T  \{ \widehat{I^{-1}(\beta_0)}  \}^{-1}  (\hat \beta - \beta_0)$
follows asymptotically $ \chi^2 (d)$, a chisquare distribution of $d$ degrees of freedom.
For $\alpha \in (0, 1)$,  let $\chi^2_{1-\alpha}(d) $ be the $(1-\alpha)$ quantile of $ \chi^2 (d)$.
A  $(1-\alpha)$-level confidence region    for $\beta_{0}$  can be constructed as
\ba
\label{CR}
\{ \beta:  n (\hat \beta - \beta)^\T \{  \widehat{I^{-1}(\beta_0)}  \}^{-1} (\hat \beta - \beta)   \leq \chi^2_{1-\alpha}(d) \}.
\ea
And
for the hypothesis  $H_0: \beta = \beta_0
\leftrightarrow H_1: \beta\neq \beta_0$,
we propose to
reject   $H_0$ at the significance level $\alpha$  if
\ba
\label{test}
n (\hat \beta - \beta_0)^\T \{  \widehat{I^{-1}(\beta_0)}  \}^{-1} (\hat \beta - \beta_0)   > \chi^2_{1-\alpha}(d).
\ea
By Theorems \ref{thm-asym-normal-linear}  and \ref{thm-variance-consistency},
the confidence region \eqref{CR}
has an asymptotically correct $(1-\alpha)$ coverage probability,
and the test defined by the rejection region \eqref{test}
has an asymptotically correct type I eror $ \alpha $.

\section{Simulations}
\label{sec-sim}
In this section, we  conduct   simulations  to assess the finite-sample performance
 of the  proposed SMPLE $\hat{\beta}$ and
the proposed confidence region \eqref{CR}  for  the linear covariate effect $\beta$.
To generate data,  we take $X$ and $Z$ to be two scalar random variables,
which are iid from the  standard normal distribution, and
take the  conditional distribution of  $T $
given  $(X, Z)$ to be an exponential distribution with mean $1/\exp(\beta_{0}X+g_{0}(Z))$.
 Therefore the conditional hazard function of $T $
given  $(X, Z)$ is  $ \exp(\beta_{0}X+g_{0}(Z))$.
We set the censoring time $C$ to follow  a uniform distribution on $(0, c)$.
We consider three scenarios   for $g_{0}$:
(I) $g_{0}(z) = -2z$,  (II) $g_{0}(z) = -|z|^{3}/2$ and  (III) $g_{0}(z) = 2|z|$.
We set $\beta_{0}  = -2$ and consider two choices for  $c$:  $5$ and $10$,
and  two sample sizes: $600$ and $800$.
A larger  $c$   results in a smaller censoring proportion.

When implementing our SMPLE,  we set  $\boldsymbol{k}_{0}$ to
be $3, 4, 3$ in Scenarios I--III, respectively.
For comparison, we also consider    the traditional Cox regression estimator  (TCR)
of $\beta$ and the partial likelihood estimator of  $\beta$
with the   $r$-order polynomial splines
 under the partially linear additive model   \citep[SPLA-$r$]{Huang1999},
 where   $r$ may be $2, 3$ or $4$.
We generate 1000 random samples to evaluate
the performance of the above five  estimation methods.

\subsection{Point estimation}
 Table \ref{tab-sim-1-0.35} presents $100$ times the simulated root mean square errors (RMSEs)
 and the absolute biases (BIASs) of these estimators.
The model assumptions of SMPLE and the SPLA-$r$ are correct in all the three scenarios,
whereas the   standard Cox model is correctly specified only in Scenario I.
As expected,   in Scenario I,  TCR has uniformly the best performance
among the five estimators under comparison in terms of RMSE and BIAS.
Nevertheless,  the SMPLE and the SPLA-$r$ estimators
have almost the same   RMSEs  and BIASs.
When the standard Cox model is misspecified in Scenarios II and III,
TCR has much larger   RMSEs and BIASs   than  SMPLE and the SPLA-$r$,
or equivalently,  SMPLE and the SPLA-$r$ have clear priority over TCR.
Compared with  SPLA-$r$,
  SMPLE is comparable and slightly inferior   in Scenarios I and II,
but  has  uniformly much smaller  RMSEs  and BIASs in Scenario III.
A possible explanation for this phenomenon is that
although   continuous in all three scenarios,
the hazard function is smooth in  Scenarios I and II  but nonsmooth in Scenario III. As the sample size $n$  or the constant $c$ increases, we have more completely observed data,
consequently  all estimators  have improved performance
when the underlying model assumption is correct.
A counterexample is  the performance of TCR in Scenarios II and III,
where TCR has larger RMSEs and BIASs as $n$ or $c$ increases.

\begin{figure}[hb]
	\begin{center}
		\centering
		\includegraphics[width=0.8\linewidth]{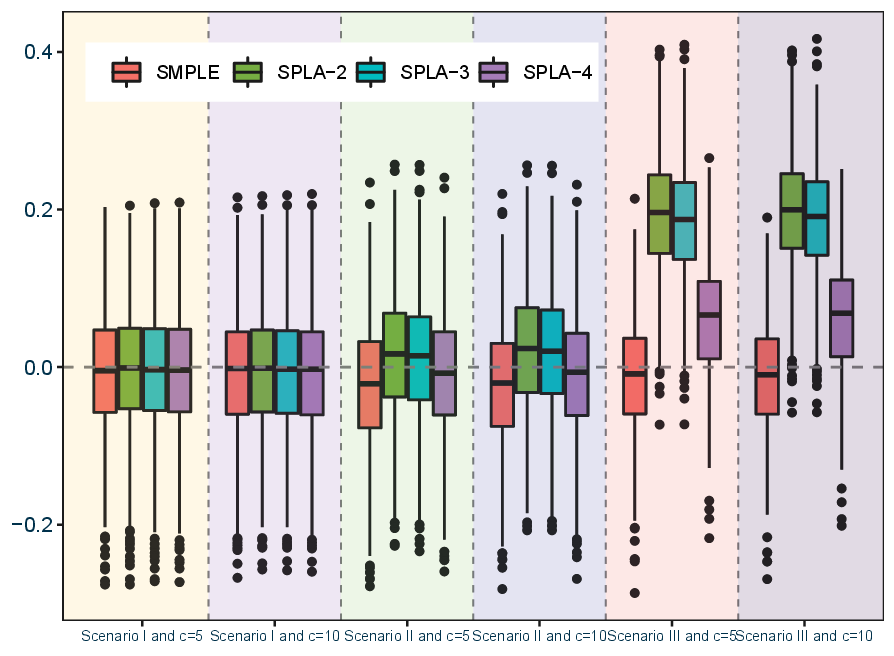}				
		\caption{Boxplots of  the  SMPLE,  SPLA-$2$,  SPLA-$3$ and SPLA-$4$  estimates (minus $\beta_0$) of $\beta$
		 when the sample size is $n=800$.}
		\label{sim-box-plot}
	\end{center}
\end{figure}

Figure \ref{sim-box-plot} displays the boxplots of
the SMPLE and SPLA-$r$  (r=2, 3, 4) estimators   (minus $\beta_0$)
of $\beta$ under study when the sample size is $n=800$.
TCR is excluded here as it has extremely large RMSEs and BIASs in Scenarios II and III.
 SMPLE and the three SPLA-$r$  exhibit almost the same performance in Scenario I.
In Scenario II, where the true hazard function is smooth,
the four methods have close variances,
but from SMPLE, to SPLA-$2$, SPLA-$3$, and SPLA-$4$,
their BIASs become smaller and smaller.
In Scenario III,  where the true hazard function is nonsmooth,
the four methods  have close variances again,
however  the three  SPLA-$r$ estimators have much larger BIASs
than SMPLE,  whose BIASs are negligible.

\begin{table}[t]
	\centering
	\caption{Simulated root mean square errors and absolute biases
 (in parentheses) of the five point estimators under comparison.
 All results have been multiplied by 100.}

\vspace{1ex}
	\begin{tabular}{ccrrrrrr}
		\toprule
		\multicolumn{1}{l}{Scenario} & \multicolumn{1}{l}{$c$} & \multicolumn{1}{l}{$n$} & \multicolumn{1}{l}{SMPLE} & \multicolumn{1}{l}{TCR} & \multicolumn{1}{l}{SPLA-$2$} & \multicolumn{1}{l}{SPLA-$3$} & \multicolumn{1}{l}{SPLA-$4$} \\
		\midrule
		\multirow{4}[2]{*}{I} & \multirow{2}[1]{*}{5} & 600   &  $\underset{(7.29)}{9.25}$    &    $\underset{(7.17)}{9.13}$ &    $\underset{(7.21)}{9.17}$   &  $\underset{(7.27)}{9.22}$     & $\underset{(7.32)}{9.29}$\\
		&       & 800  &   $\underset{(6.36)}{8.04}$   & $\underset{(6.30)}{7.95}$    &    $\underset{(6.33)}{7.99}$  &  $\underset{(6.35)}{8.01}$      & $\underset{(6.36)}{8.03}$ \\
		& \multirow{2}[1]{*}{10} & 600  &    $\underset{(7.19)}{8.96}$   &  $\underset{(7.07)}{8.86}$     &   $\underset{(7.11)}{8.90}$    & $\underset{(7.17)}{8.96}$    & $\underset{(7.21)}{9.01}$ \\
		&       & 800  &   $\underset{(6.07)}{7.65}$   & $\underset{(6.03)}{7.57}$     &  $\underset{(6.05)}{7.61}$    &     $\underset{(6.05)}{7.62}$  & $\underset{(6.06)}{7.64}$\\
		\midrule
		\multirow{4}[2]{*}{II} & \multirow{2}[1]{*}{5} & 600   &  $\underset{(8.19)}{10.29}$      &   $\underset{(67.57)}{68.89}$     &    $\underset{(7.65)}{9.46}$    &  $\underset{(7.66)}{9.46}$      &  $\underset{(7.78)}{9.70}$ \\
		&       & 800   &   $\underset{(6.68)}{8.49}$    &   $\underset{(69.02)}{70.06}$    &   $\underset{(6.53)}{8.17}$    &   $\underset{(6.44)}{8.10}$    & $\underset{(6.36)}{8.09}$ \\
		&  \multirow{2}[1]{*}{10} & 600   &    $\underset{(7.86)}{9.89}$    &     $\underset{(77.53)}{78.58}$   &     $\underset{(7.42)}{9.21}$   &     $\underset{(7.40)}{9.18}$   &  $\underset{(7.44)}{9.30}$ \\
		&       & 800   &   $\underset{(6.46)}{8.20}$    &    $\underset{(78.79)}{79.62}$   &    $\underset{(6.48)}{8.04}$   &     $\underset{(6.38)}{7.93}$  & $\underset{(6.15)}{7.78}$ \\
		\midrule
		\multirow{4}[2]{*}{III} & \multirow{2}[1]{*}{5} & 600  & $\underset{(6.87)}{8.54}$      &  $\underset{(63.15)}{63.46}$     & $\underset{(18.83)}{20.54}$      &   $\underset{(17.63)}{19.39}$    & $\underset{(7.95)}{9.74}$ \\
		&       & 800   &  $\underset{(5.74)}{7.28}$     &   $\underset{(63.35)}{63.61}$    &     $\underset{(19.60)}{20.92}$  &   $\underset{(18.62)}{19.98}$    & $\underset{(7.74)}{9.28}$ \\
		& \multirow{2}[1]{*}{10} & 600   &   $\underset{(6.78)}{8.43}$    &   $\underset{(63.40)}{63.70}$    &   $\underset{(19.40)}{21.05}$    &   $\underset{(18.21)}{19.91}$    & $\underset{(8.04)}{9.86}$ \\
		&       & 800   &   $\underset{(5.66)}{7.17}$    &  $\underset{(63.47)}{63.71}$     &    $\underset{(20.01)}{21.25}$   &     $\underset{(19.04)}{20.31}$  &  $\underset{(7.78)}{9.34}$\\
		\bottomrule
	\end{tabular}%
	\label{tab-sim-1-0.35}%
\end{table}%

\begin{table}[ht]
	\centering
	\caption{Simulated coverage probabilities and average interval lengths (in parentheses)
		of the confidence intervals at the 95\% confidence level based on the five
		estimators under comparison}
	\begin{tabular}{cccccccc}
		\toprule
		\multicolumn{1}{l}{Scenario} & \multicolumn{1}{l}{$c$} & \multicolumn{1}{l}{$n$} & \multicolumn{1}{l}{SMPLE} & \multicolumn{1}{l}{TCR} & \multicolumn{1}{l}{SPLA-$2$} & \multicolumn{1}{l}{SPLA-$3$} & \multicolumn{1}{l}{SPLA-$4$} \\
		\midrule
		\multirow{4}[2]{*}{I} & \multirow{2}[1]{*}{5} & 600   &   $\underset{(0.412)}{0.947}$   &  $\underset{(0.384)}{0.930}$   &   $\underset{(0.396)}{0.940}$    &  $\underset{(0.409)}{0.944}$   &  $\underset{(0.422)}{0.948}$\\
		&       & 800   &   $\underset{(0.343)}{0.941}$   &   $\underset{(0.324)}{0.937}$    &  $\underset{(0.332)}{0.937}$      &   $\underset{(0.340)}{0.941}$    &  $\underset{(0.348)}{0.949}$\\
		& \multirow{2}[1]{*}{10} & 600   &  $\underset{(0.388)}{0.940}$      &   $\underset{(0.364)}{0.930}$   &   $\underset{(0.374)}{0.937}$     &  $\underset{(0.385)}{0.940}$   &  $\underset{(0.395)}{0.942}$  \\
		&       & 800  &  $\underset{(0.327)}{0.945}$     &  $\underset{(0.310)}{0.935}$     &  $\underset{(0.317)}{0.941}$    &   $\underset{(0.325)}{0.949}$    & $\underset{(0.332)}{0.948}$ \\
		\midrule
		\multirow{4}[2]{*}{II} & \multirow{2}[1]{*}{5} & 600  &    $\underset{(0.459)}{0.947}$    &    $\underset{(0.452)}{0.001}$    &    $\underset{(0.410)}{0.941}$    &   $\underset{(0.422)}{0.951}$     &  $\underset{(0.440)}{0.951}$ \\
		&       & 800   &   $\underset{(0.382)}{0.950}$     &   $\underset{(0.396)}{0}$     &   $\underset{(0.346)}{0.926}$     &   $\underset{(0.355)}{0.938}$     & $\underset{(0.368)}{0.949}$  \\
		& \multirow{2}[1]{*}{10} & 600   &    $\underset{(0.429)}{0.943}$    &  $\underset{(0.437)}{0}$      &    $\underset{(0.383)}{0.929}$    &   $\underset{(0.393)}{0.935}$    &  $\underset{(0.407)}{0.938}$ \\
		&       & 800   &   $\underset{(0.355)}{0.945}$    &   $\underset{(0.383)}{0}$    &     $\underset{(0.321)}{0.919}$  &   $\underset{(0.329)}{0.931}$    & $\underset{(0.340)}{0.942}$  \\
		\midrule
		\multirow{4}[2]{*}{III} & \multirow{2}[1]{*}{5} & 600  &    $\underset{(0.384)}{0.955}$   &  $\underset{(0.261)}{0}$     &    $\underset{(0.339)}{0.424}$   &   $\underset{(0.350)}{0.494}$    & $\underset{(0.364)}{0.908}$ \\
		&       & 800   &   $\underset{(0.321)}{0.955}$    &    $\underset{(0.223)}{0}$   &   $\underset{(0.289)}{0.260}$    &    $\underset{(0.296)}{0.302}$   & $\underset{(0.304)}{0.863}$ \\
		& \multirow{2}[1]{*}{10} & 600   &   $\underset{(0.373)}{0.952}$    &  $\underset{(0.254)}{0}$     &  $\underset{(0.330)}{0.370}$     &    $\underset{(0.340)}{0.441}$   & $\underset{(0.353)}{0.893}$ \\
		&       & 800   &   $\underset{(0.312)}{0.948}$    &    $\underset{(0.217)}{0}$   &   $\underset{(0.281)}{0.233}$    &   $\underset{(0.288)}{0.279}$    & $\underset{(0.296)}{0.847}$ \\
		\bottomrule
	\end{tabular}%
	\label{tab-sim-2-0.35}%
\end{table}%

\subsection{Interval estimation}

We continue to compare  the proposed confidence interval (CI) for $\beta$
with those based on TCR and the three SPLA-$r$ estimators.
To be specific, let $\tilde \beta$ denote a generic
  estimator  of $\beta$  and
let $\sigma^2$ denote its asymptotic variance.
We use the proposed data-splitting method  with $\tilde{\alpha} = 0.35$ to estimate
the asymptotic variance $\sigma^2$, and let $\tilde \sigma^2$ be the resulting estimator.
Then a Wald-type confidence interval at the confidence level 95\% for $\beta$ is
$\tilde \beta \pm  n^{-1/2} \tilde \sigma u_{0.975}$
and its length is $ 2 n^{-1/2} \tilde \sigma u_{0.975}$,
where $u_{0.975}$ is the $0.975$ quantile of the standard normal distribution.
 Table \ref{tab-sim-2-0.35} presents the simulated coverage probabilities and
average interval lengths  of the Wald-type confidence intervals
based on the five estimators.

\begin{figure}[ht]
	\begin{center}
		\centering
				\includegraphics[width=0.98\linewidth]{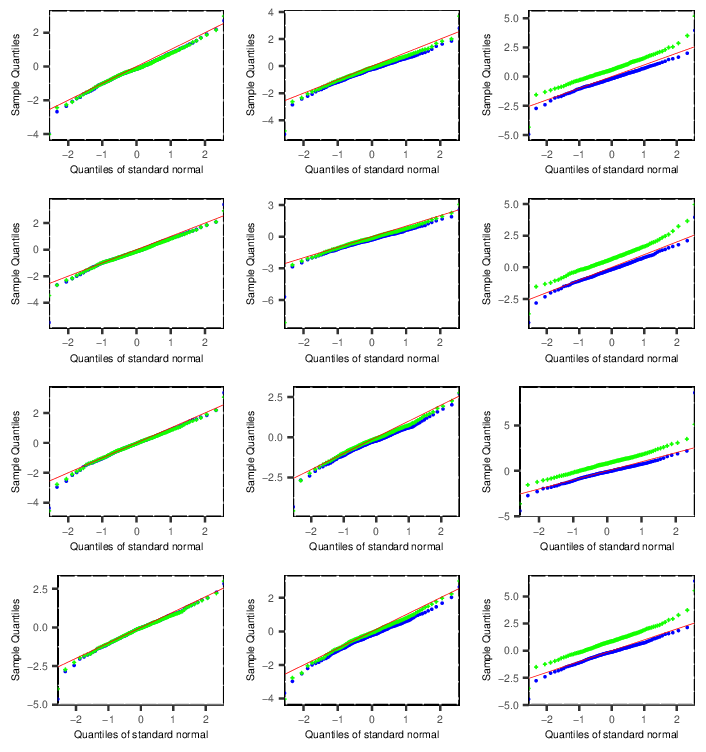}				
		\caption{QQ-plots of    the standardized estimates
($
 \sqrt{n}\{ \widehat{I^{-1}(\beta_{0})}   \}^{-1/2} (\tilde{\beta}-\beta_{0} )
$)
of   SMPLE (blue circle) and SPLA-$4$ (green cross) versus $N(0,1)$.
Columns 1--3: Scenarios I--III;  Row 1: $c=5$ and $n=600$;
 Row 2: $c=10$ and $n=600$;  Row 3: $c=5$ and $n=800$;  Row 4: $c=10$ and $n=800$. }
		\label{sim-qq-plot}
	\end{center}
\end{figure}

We observe that  in all the three scenarios,
the SMPLE-based CI  always has very accurate
coverage probabilities.  In contrast,
the TCR-based CI has desirable coverage accuracy
only in Scenario I, where the Cox model is correct,
and its performance in Scenarios II and III
is   unacceptable.
The CIs based on SPLA-$r$ have acceptable
coverage accuracy in Scenarios I and II,
although they may have slight undercoverage,
but performs very poor in Scenario III.
Interestingly,  the  CI  based on SPLA-$r$
with a larger $r$ have more accurate coverage accuracy.
To get more insights about these simulation results,
we display in Figure \ref{sim-qq-plot}
the QQ-plots  of  the standardized estimates
$
 \sqrt{n}\{ \widehat{I^{-1}(\beta_{0})}   \}^{-1/2} (\tilde{\beta}-\beta_{0} )
$
for $\tilde{\beta}$ being the SMPLE and SPLA-$4$ estimators.
We see that the distribution of the standardized  SMPLE estimator
is always very close to the standard normal.
In contrast, the distribution of the standardized  SPLA-$4$  estimator
is  close to the standard normal in Scenarios I and II,
but far away from the standard normal in Scenario  III.
This may explains why the SMPLE-based CI always has
desirable coverage accuracy but the  SPLA-$4$-based  CI
has severely undercoverage in Scenario III.

In summary, when the Cox model is correct,
 the  point and interval estimators  based on the proposed  SMPLE
have very close performance as those based on TCR,
which is  the optimal estimation method under the Cox model.
And they  have obvious advantages over  those based on TCR
or SPLA-$r$  when the Cox model is incorrect or the hazard function is nonsmooth.
Another advantage of SMPLE over SPLA-$r$  is that
it is free from tuning parameters but the latter is not.
Tuning parameters  usually have a big influence on
the subsequent analysis but their optimal choices are difficult
to determine in practice.

%

\section{A real application}
\label{sec-real-data}
We apply the proposed SMPLE method to
analyze the Rotterdam Breast Cancer (RBC) dataset,
which is publicly available from the R package {\tt survival}.
This dataset was used in \cite{royston2013external} to
perform external validation of a Cox prognostic model.
The RBC dataset comprises  2982 primary breast cancers patients
whose records   were included in the Rotterdam tumor bank.
We focus on studying potential factors that affect the
survival time ($T$) of breast cancers patients,
defined as the days from primary surgery to the death.
We consider four covariates: progesterone receptors (unit: fmol/l, $X$),
age at surgery ($Z_{1}$), number of positive lymph nodes ($Z_{2}$),
estrogen receptors (unit: fmol/l, $Z_{3}$).  Let $Z=(Z_1, Z_2, Z_3)^\T$. All covariates have been scaled to between $0$ and $1$. The impact of the content  of progesterone receptors
 on the treatment of breast cancer is not clear.
We consider modelling the conditional hazard of $T$ given $(X=x, Z=z)$
by the following partially linear Cox additive model
\ba
\label{real-model}
\lambda_{T}(t\mid x, z) = \lambda(t)\exp\{ x \beta+f_{1}(z_{1})+f_{2}(z_{2})+f_{3}(z_{3}) \},
\ea
where $\lambda(\cdot)$, $f_{1}(\cdot)$  $f_{2}(\cdot)$  and  $f_{3}(\cdot)$
are all left un-specified. Since elderly patients are typically
predisposed to a higher risk of cancer deterioration or recurrence,
we postulate that $f_{1}(\cdot)$  is  monotonically increasing.
Similarly,  $f_{2}(\cdot)$ and $f_{3}(\cdot)$ are  also assumed to be monotonically increasing functions. Inspired by our extensive simulation results and discussion in Remark \ref{remark-convex-monotone}, we further impose
a  convexity or concavity restriction on $f_1(\cdot)$--$f_3(\cdot)$.

Because menopause generally affects the therapeutic effect and
the mortality of breast cancer,
we   categorize the patients into two groups
based on whether they have reached menopause or not.
 Group I (premenopausal) comprises $1312$ samples,
 exhibiting a censor rate of $64.33\%$,
 while Group II (postmenopausal) consists of
 $1670$ samples with a censor rate of $51.86\%$.
 We   apply the point and interval estimation methods considered in the simulation section
to estimate the linear covariate effect under model \eqref{real-model}.
We set the confidence level   to  95\% when constructing confidence intervals.
In the proposed variance estimation method, we take
 $\tilde{\alpha} = 0.3$  and construct variance estimates by
  repeating the splitting and estimating procedure twenty times
  and taking average.
  Table  \ref{tab-real-data} summarizes the analysis results.

\begin{table}[ht]
	\caption{
		Results of point and interval estimation for the RBC dataset. 	\label{tab-real-data} }
	\tabcolsep 2pt
	\centering
	\begin{tabular}{cccccc}
		\toprule
		& \multicolumn{2}{c}{Group I} & & \multicolumn{2}{c}{Group II} \\
		& Point estimate & Interval estimate & & Point estimate & Interval estimate \\
		\cline{2-3} \cline{5-6}
		SMPLE & $-5.96$ & $[-8.96, -2.95]$  & & $-0.62$ & $[-1.23, -0.01]$    \\
		TCR & $-5.19$ & $[-8.22, -2.16]$ & & $-0.70$ & $[-1.41, 0.01]$  \\
		SPLA-$2$ & $-5.33$ & $[-8.48, -2.18]$ &  &$-0.39$ & $[-1.08, 0.30]$   \\
		SPLA-$3$ & $-4.98$ & $[-8.26, -1.69]$  &  &$-0.40$ &$[-1.13, 0.33]$  \\
		SPLA-$4$ & $-5.15$ & $[-8.66, -1.64]$ & &$-0.34$ & $[-1.08, 0.39]$  \\
		\bottomrule
	\end{tabular}
\end{table}

The point estimates of $\beta$ are about $-5$ or less  in Group I,
which are much less than  those in Group II, namely $-0.7$ or larger.
This implies that the effect of progesterone receptors ($X$)
on the hazard rate of survival time among breast cancer patients
has a significant increase  in Group II   than in Group I,
which potentially stems from the substantial reduction in progesterone levels
in women's bodies post-menopause.
For patients who have not reached menopause, the $95\%$ confidence intervals derived
from all five methods exclude $0$,
suggesting a significant effect of the progesterone receptors
on the survival time of breast cancer patients.
Among these techniques, our proposed SMPLE method
produces the shortest confidence interval.
Conversely, for postmenopausal patients,
only  our SMPLE-based confidence interval  excludes $0$,
which together with the point estimate  $-0.62$
suggests that the progesterone receptors still has
a significant effect in  increasing the survival time  of breast cancer patients.
Once again, our confidence interval exhibits the shortest confidence interval. According to our reasonable shape restriction assumptions and our simulation experience especially those in Scenario III,
we believe that the analysis results based on our SMPLE are
more reliable, which  also demonstrate its priority over
the competing methods.

\begin{figure}[ht]
	\begin{center}
		\centering
		\includegraphics[width=1\linewidth]{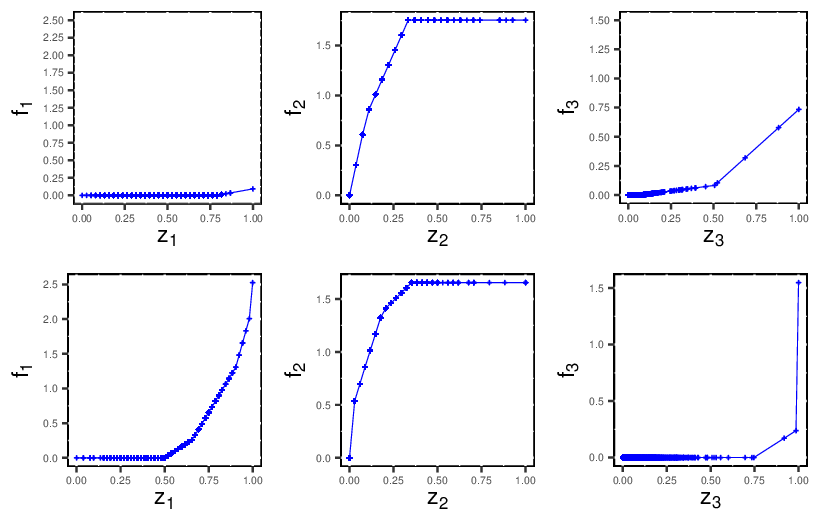}				
		\caption{Curves of the fitted $f_{1}(\cdot)$--$f_{3}(\cdot)$ in model \eqref{real-model}
by the proposed shape-restricted maximum likelihood estimation. Row 1: Group I, Row 2: Group II. }
		\label{real-plot}
	\end{center}
\end{figure}

Figure \ref{real-plot} depicts the estimates of $f_{1}(\cdot)$--$f_{3}(\cdot)$
in model \eqref{real-plot} under our shape constraints.
All the estimated curves are continuous and piecewise linear   functions
in both groups.
We see that  both number of positive lymph nodes ($Z_2$) and  estrogen receptors ($Z_3$)  have
very similar effects   in the two groups on the survival time $T$.
Exceptionally,  the covariate age at the time of surgery ($Z_1$)
  has a significantly more pronounced effect on $T$
  in Group II than in Group I.
   This discrepancy may stem from the fact that premenopausal women
   tend to be younger, and thus, age has not yet exerted a substantial influence on their physical condition.


\section{Concluding remarks}
\label{sec-con}
In this paper, we study an additive and shape-restricted partially linear Cox  model.
We systematically investigate the consistencies and convergence rates
of the SMPLEs for  the additive components and the baseline
cumulative hazard function.
 Notably, the convergence rates of the SMPLEs  for the infinite-dimensional
  parameters, including the additive components and $\Lambda(\cdot)$,
are independent of the covariate dimension, exhibiting an elegant form.
We find that the SMPLE for $\beta$ is $\sqrt{n}$ consistent
and asymptotically semiparametric efficiency.

To estimate the asymptotic variance, we propose a flexible variance estimation method
by data splitting.
We prove that the resulting variance estimator is consistent
if $\tilde{\alpha} \in (0, 1)$, namely,
both the number of split  samples and the sample size of each split sample
all goes to infinity.
 With this variance estimator, valid confidence regions for
the linear covariate effect can be constructed.
We set $\tilde{\alpha}$ to 0.3 around in our numerical studies.
Natural questions are if there is an optimal choice of $\tilde{\alpha}$,
and if so, how to determine the optimal value.
For the time being, we do not have any results on these questions,
which may be left  for future research.

In addition to the  efficiency studies of the
estimation of  the linear covariate effect,
there is also a compelling interest in developing the minimax lower bound
for the estimation of the additive  and shape-restricted
 components.
We conjecture that, by leveraging the techniques
in \citep{guntuboyina2015global, Zhong2022},
the minimax lower bound under the $L_{2}(P)$ norm
for the convex/concave component
may be  (asymptotically) lower bounded by $n^{-2/5}$
up to a constant factor.
Furthermore, if none of additive components is monotonic,
the  SMPLEs for the convex/concave components are rate-optimal,
as suggested by Theorem \ref{thm-rate-convergence-SMPLE}.
However, a formal proof of this result is beyond the scope of this paper.

\bibliographystyle{natbib}
\bibliography{cox-additive}

\begin{thebibliography}{}

\bibitem[Balabdaoui {\em et~al.}(2019)Balabdaoui, Groeneboom, and
  Hendrickx]{balabdaoui2019score}
Balabdaoui, F., Groeneboom, P., and Hendrickx, K. (2019).
\newblock Score estimation in the monotone single-index model.
\newblock {\em Scandinavian Journal of Statistics\/}, {\bf 46}(2), 517--544.

\bibitem[Bickel {\em et~al.}(1993)Bickel, Klaassen, Ritov, and
  Wellner]{bickel1993efficient}
Bickel, P.~J., Klaassen, J., Ritov, Y., and Wellner, J.~A. (1993).
\newblock {\em Efficient and adaptive estimation for semiparametric models\/},
  volume~4.
\newblock Springer.

\bibitem[Chang {\em et~al.}(2007)Chang, Chien, Hsiung, Wen, and
  Wu]{chang2007shape}
Chang, I.-S., Chien, L.-C., Hsiung, C.~A., Wen, C.-C., and Wu, Y.-J. (2007).
\newblock Shape restricted regression with random bernstein polynomials.
\newblock {\em Lecture Notes-Monograph Series\/}, pages 187--202.

\bibitem[Chen and Samworth(2016)Chen and Samworth]{ChenandSamworth2016}
Chen, Y. and Samworth, R.~J. (2016).
\newblock Generalized additive and index models with shape constraints.
\newblock {\em Journal of the Royal Statistical Society\/}, {\bf 78}(4),
  729--754.

\bibitem[Cheng(2009)Cheng]{Cheng2009}
Cheng, G. (2009).
\newblock Semiparametric additive isotonic regression.
\newblock {\em Journal of Statistical Planning and Inference\/}, {\bf 139}(6),
  1980--1991.

\bibitem[Cox(1972)Cox]{Cox1972}
Cox, D.~R. (1972).
\newblock Regression models and life-tables.
\newblock {\em Journal of the Royal Statistical Society: Series B
  (Methodological)\/}, {\bf 34}(2), 187--202.

\bibitem[Cox(1975)Cox]{Cox1975}
Cox, D.~R. (1975).
\newblock Partial likelihood.
\newblock {\em Biometrika\/}, {\bf 62}(2), 269--276.

\bibitem[Deng {\em et~al.}(2023)Deng, Xu, Fu, Wang, and Qin]{deng2023active}
Deng, G., Xu, G., Fu, Q., Wang, X., and Qin, J. (2023).
\newblock Active-set algorithm-based statistical inference for shape-restricted
  generalized additive cox regression models.
\newblock {\em Journal of Statistical Computation and Simulation\/}, {\bf
  93}(3), 416--441.

\bibitem[Deng and Zhang(2020)Deng and Zhang]{DengandZhang2020}
Deng, H. and Zhang, C.-H. (2020).
\newblock {Isotonic regression in multi-dimensional spaces and graphs}.
\newblock {\em The Annals of Statistics\/}, {\bf 48}(6), 3672--3698.

\bibitem[Feng {\em et~al.}(2022)Feng, Chen, Han, Carroll, and
  Samworth]{feng2022nonparametric}
Feng, O.~Y., Chen, Y., Han, Q., Carroll, R.~J., and Samworth, R.~J. (2022).
\newblock Nonparametric, tuning-free estimation of s-shaped functions.
\newblock {\em Journal of the Royal Statistical Society Series B: Statistical
  Methodology\/}, {\bf 84}(4), 1324--1352.

\bibitem[Groeneboom and Hendrickx(2017)Groeneboom and
  Hendrickx]{Groeneboom2017}
Groeneboom, P. and Hendrickx, K. (2017).
\newblock {The nonparametric bootstrap for the current status model}.
\newblock {\em Electronic Journal of Statistics\/}, {\bf 11}(2), 3446--3484.

\bibitem[Groeneboom and Jongbloed(2014)Groeneboom and
  Jongbloed]{groeneboom2014}
Groeneboom, P. and Jongbloed, G. (2014).
\newblock {\em Nonparametric Estimation under Shape Constraints: Estimators,
  Algorithms and Asymptotics\/}.
\newblock London, UK: Cambridge University Press.

\bibitem[Guntuboyina and Sen(2015)Guntuboyina and Sen]{guntuboyina2015global}
Guntuboyina, A. and Sen, B. (2015).
\newblock Global risk bounds and adaptation in univariate convex regression.
\newblock {\em Probability Theory and Related Fields\/}, {\bf 163}(1),
  379--411.

\bibitem[Hastie and Tibshirani(1986)Hastie and Tibshirani]{1986Generalized}
Hastie, T.~J. and Tibshirani, R.~J. (1986).
\newblock Generalized additive models.
\newblock {\em Statistical Science\/}, {\bf 1}(3), 297--310.

\bibitem[Hastie and Tibshirani(1990)Hastie and Tibshirani]{1990Generalized}
Hastie, T.~J. and Tibshirani, R.~J. (1990).
\newblock {\em Generalized Aditive Models\/}.
\newblock London: Chapman and Hall.

\bibitem[Heller(2001)Heller]{heller2001cox}
Heller, G. (2001).
\newblock The cox proportional hazards model with a partly linear relative risk
  function.
\newblock {\em Lifetime data analysis\/}, {\bf 7}, 255--277.

\bibitem[Horowitz and Lee(2017)Horowitz and Lee]{horowitz2017nonparametric}
Horowitz, J.~L. and Lee, S. (2017).
\newblock Nonparametric estimation and inference under shape restrictions.
\newblock {\em Journal of Econometrics\/}, {\bf 201}(1), 108--126.

\bibitem[Huang(1999)Huang]{Huang1999}
Huang, J. (1999).
\newblock {Efficient estimation of the partly linear additive Cox model}.
\newblock {\em The Annals of Statistics\/}, {\bf 27}(5), 1536--1563.

\bibitem[Huang(2002)Huang]{Huang2002}
Huang, J. (2002).
\newblock A note on estimating a partly linear model under monotonicity
  constraints.
\newblock {\em Journal of Statistical Planning and Inference\/}, {\bf
  107}(1-2), 343--351.

\bibitem[Kuchibhotla {\em et~al.}(2023)Kuchibhotla, Patra, and
  Sen]{Kuchibhotla2023}
Kuchibhotla, A.~K., Patra, R.~K., and Sen, B. (2023).
\newblock Semiparametric efficiency in convexity constrained single index
  model.
\newblock {\em Journal of the American Statistical Association\/}, {\bf
  118}(541), 272--286.

\bibitem[Matzkin(1991)Matzkin]{matzkin1991semiparametric}
Matzkin, R.~L. (1991).
\newblock Semiparametric estimation of monotone and concave utility functions
  for polychotomous choice models.
\newblock {\em Econometrica\/}, {\bf 59}(5), 1315--1327.

\bibitem[O'Sullivan(1993)O'Sullivan]{o1993nonparametric}
O'Sullivan, F. (1993).
\newblock Nonparametric estimation in the cox model.
\newblock {\em The Annals of Statistics\/}, {\bf 21}(1), 124--145.

\bibitem[Qin {\em et~al.}(2014)Qin, Garcia, Ma, Tang, Marder, and
  Wang]{qin2014combining}
Qin, J., Garcia, T.~P., Ma, Y., Tang, M.-X., Marder, K., and Wang, Y. (2014).
\newblock Combining isotonic regression and em algorithm to predict genetic
  risk under monotonicity constraint.
\newblock {\em The annals of applied statistics\/}, {\bf 8}(2), 1182--1208.

\bibitem[Qin {\em et~al.}(2021)Qin, Deng, Ning, Yuan, and
  Shen]{qin2021estrogen}
Qin, J., Deng, G., Ning, J., Yuan, A., and Shen, Y. (2021).
\newblock Estrogen receptor expression on breast cancer patients' survival
  under shape-restricted cox regression model.
\newblock {\em The annals of applied statistics\/}, {\bf 15}(3), 1291--1307.

\bibitem[Rong {\em et~al.}(2024)Rong, Zhao, Zheng, and Li]{rong2024kernel}
Rong, Y., Zhao, S.~D., Zheng, X., and Li, Y. (2024).
\newblock Kernel cox partially linear regression: Building predictive models
  for cancer patients' survival.
\newblock {\em Statistics in Medicine\/}, {\bf 43}(1), 1--15.

\bibitem[Royston and Altman(2013)Royston and Altman]{royston2013external}
Royston, P. and Altman, D.~G. (2013).
\newblock External validation of a cox prognostic model: principles and
  methods.
\newblock {\em BMC medical research methodology\/}, {\bf 13}, 1--15.

\bibitem[Sasieni(1992)Sasieni]{Sasieni1992}
Sasieni, P. (1992).
\newblock Information bounds for the conditional hazard ratio in a nested
  family of regression models.
\newblock {\em Journal of the Royal Statistical Society: Series B
  (Methodological)\/}, {\bf 54}(2), 617--635.

\bibitem[Sleeper and Harrington(1990)Sleeper and
  Harrington]{sleeper1990regression}
Sleeper, L.~A. and Harrington, D.~P. (1990).
\newblock Regression splines in the cox model with application to covariate
  effects in liver disease.
\newblock {\em Journal of the American Statistical Association\/}, {\bf
  85}(412), 941--949.

\bibitem[Varian(1984)Varian]{varian1984nonparametric}
Varian, H.~R. (1984).
\newblock The nonparametric approach to production analysis.
\newblock {\em Econometrica\/}, {\bf 52}(3), 579--597.

\bibitem[Wang and Ghosh(2012)Wang and Ghosh]{wang2012shape}
Wang, J. and Ghosh, S.~K. (2012).
\newblock Shape restricted nonparametric regression with bernstein polynomials.
\newblock {\em Computational Statistics \& Data Analysis\/}, {\bf 56}(9),
  2729--2741.

\bibitem[Zhong {\em et~al.}(2022)Zhong, Mueller, and Wang]{Zhong2022}
Zhong, Q., Mueller, J., and Wang, J.-L. (2022).
\newblock {Deep learning for the partially linear Cox model}.
\newblock {\em The Annals of Statistics\/}, {\bf 50}(3), 1348--1375.

\end{thebibliography}

\end{document}